\title{\bf{A new infinite family of irregular algebraic surfaces with canonical map of degree 8}}
\author{
	NGUYEN BIN\\
}
\date{\today}

\newcommand{\Addresses}{{% additional braces for segregating \footnotesize
		\bigskip
		\footnotesize
		    \text{Mathematics Division,}\par\nopagebreak
			\text{National Center for Theoretical Sciences,}\par\nopagebreak	
			\text{National Taiwan University,}\par\nopagebreak	
			\text{Taiwan.}\par\nopagebreak		
		\textit{E-mail address}: \texttt{nguyenbin@ncts.ntu.edu.tw}
				
	}}

\newcommand\blfootnote[1]{%
	\begingroup
	\renewcommand\thefootnote{}\footnote{#1}%
	\addtocounter{footnote}{-1}%
	\endgroup
}

\date{\today}

\documentclass[12pt]{article}
\usepackage{amsmath, amsthm, amssymb, mathrsfs, amscd, amsthm, amscd,amsfonts,calligra,mathrsfs,adjustbox,lipsum}
\usepackage{todonotes}
\usepackage{indentfirst,url,xypic}
\usepackage[all]{xy}
\usepackage{url}
\usepackage{tikz}
\usepackage[colorlinks,plainpages]{hyperref}
\hypersetup{
	colorlinks=true,
	linkcolor=blue,
	filecolor=magenta,      
	urlcolor=cyan,
}

\DeclareMathOperator{\degree}{deg}
\DeclareMathOperator{\image}{im}

\newtheorem{Theorem}{Theorem}
\newtheorem{Proposition}{Proposition }

\newtheorem{Question}{Question}

\theoremstyle{remark}
\newtheorem{Remark}{Remark}

\newcommand{\MSC}{\textbf{Mathematics Subject Classification (2010):}}

\newcommand{\Key}{\textbf{Key words:}}

\begin{document}
\maketitle
\begin{abstract}
	In this note we construct an unlimited family of irregular algebraic surfaces of general type with canonical map of degree $ 8 $, irregularity $ 1 $ and arbitrarily large geometric genus such that the image of the canonical map is not a surface of minimal degree.
\end{abstract}

\blfootnote{\MSC{ 14J29}.}
\blfootnote{\Key{ Surfaces of general type, Canonical maps, Abelian covers.}}

\section{Introduction} 
    Let $ X $ be a minimal smooth complex surface of general type and denote by $ \xymatrix{\varphi_{\left| K_X\right| }:X \ar@{.>}[r] & \mathbb{P}^{p_g\left( X\right)-1}} $ the canonical map of $ X $, where $ K_{X} $ is the canonical divisor of $ X $ and $ p_g\left( X\right) = \dim H^0\left( X, K_{X}\right)  $ is the geometric genus. The existence of surfaces of general type with non-birational canonical map has been studied intensively by many authors in the last decades. This problem is motivated by the work of A. Beauville \cite{MR553705}. We refer to the recent preprint by M. Mendes Lopes and R. Pardini \cite{lopes2021degree} on the subject. We restrict our attention to the existence of surfaces of general type with canonical map of degree $ 8 $. It is known that for surfaces of general type, the degree $ d $ of the canonical map is at most $ 9 $ if the holomorphic Euler-Poincar\'{e} characteristic $ \chi\left( \mathcal{O}_X\right)  $ is bigger than or equal to $ 31 $. This was proven by A. Beauville in {\cite{MR553705}}. In 1986, G. Xiao improved this result by showing that the degree of the canonical map is less than or equal to $ 8 $ if the geometric genus of the surface is bigger than $ 132 $ \cite{MR842626}. The first unlimited family of surfaces with canonical map of degree $ 8 $ was found by A. Beauville \cite{MR4131989}. These examples were constructed as double covers of the product surface of a non-hyperelliptic curve of genus $ 3 $ and the rational curve $ \mathbb{P}^1 $. In \cite{MR4131989}, some unlimited families of surfaces with canonical map of degree $ 8 $ were found by as $ \mathbb{Z}_{2}^3$-covers of the first Hirzebruch surface $ \mathbb{F}_1 $ or its blow-up in a point. In these all known examples, the image of the canonical map is a surface of minimal degree. In this note we construct an unlimited family of algebraic surfaces with $ d = 8 $ and arbitrarily large $ p_g $ such that the image of the canonical map is not of minimal degree.    
    
   \begin{Theorem}\label{the main theorem}
   	Let $ n  $ be an integer with $ n \ge 3 $. There exist minimal surfaces of general type $ X $ satisfying	
   	$$
   	\begin{tabular}{|c|c|c|c|c|}
   		\hline
   		$ K_X^2 $ &$ p_g\left( X\right)  $&$ q\left( X\right)  $ &$ \degree\left( \image\varphi_{\left| K_X\right| }  \right)$&$ \left| K_X \right| $ \\ \hline
   		$ 16n $ & $ 2n $ & $ 1 $ &$ 2n $& is base point free \\ \hline
   	\end{tabular}
   	$$
   	such that the canonical map $ \varphi_{\left| K_X  \right|}  $ has degree $ 8 $.
   \end{Theorem}
    
    \noindent    
    In the above theorem, $ q\left( X\right)= \dim H^1\left( X, K_{X}\right)$ is the irregularity of $ X $. These surfaces are constructed as $ \mathbb{Z}_{2}^3$-covers of the product surface of a smooth elliptic curve $ C $ and the rational curve $ \mathbb{P}^1 $. The building data $ \left\lbrace L_{\chi}, D_{\sigma} \right\rbrace_{\chi,\sigma}$ of $ \mathbb{Z}_{2}^3$-covers (see Section \ref{z_2^3-coverings}) are chosen such that there is a character $ \chi' $ of $ \mathbb{Z}_2^3 $ with arbitrarily large $ h^{0}\left( L_{\chi'} + K_{\mathbb{P}^1 \times C}\right)  $ and that  $ h^{0}\left( L_{\chi} + K_{\mathbb{P}^1 \times C}\right)  $ vanishes for all other characters $ \chi $ of $ \mathbb{Z}_2^3 $. From the decomposition of the space of $ 2 $-forms of the surfaces (see Proposition \ref{invariants of Z_2^3 cover})
    \begin{align*}
    	H^{0}\left( X, K_X\right) = H^{0}\left( \mathbb{P}^1 \times C, K_{\mathbb{P}^1 \times C}\right) \oplus \bigoplus_{\chi \ne  \chi_{000}}{H^{0}\left( \mathbb{P}^1 \times C, K_{\mathbb{P}^1 \times C} +L_{\chi}\right)}, 
    \end{align*}
    \noindent
    such a choice of the building data allows to conclude that the canonical map of $ X $ factors through the $ \mathbb{Z}_{2}^3$-cover. Furthermore, we choose the divisor $ L_{\chi'} $ in $ \mathbb{P}^1 \times C $ such that $ L_{\chi'} + K_{\mathbb{P}^1 \times C} \equiv E + \sum\limits_{i=1}^{n}{F_{ii}} $, where $ E $ is a general elliptic fiber and $ F_{ii} $ are distinct rational fibres of $ \mathbb{P}^1 \times C $. This leads to the fact that the map $ \varphi_{\left| E + \sum\limits_{i=1}^{n}{F_{ii}} \right| } $ is of degree one for all $ n \ge 3 $. Thus, the canonical map of $ X $ is of degree $ 8 $. We notice that the covers of such a product guarantees the irregularity of the result by pulling back $ 1 $-forms.
    
    In our construction, if $ n= 2 $, we obtain a smooth minimal surface of general type $ X $ with $ K_X^2 = 32, p_g\left( X\right) =4, q\left( X\right) =1 $ and $ d = 16 $ since the linear system $ \left| E + F_{11} + F_{22} \right|  $ is a map of degree $ 2 $ onto $ \mathbb{P}^1\times \mathbb{P}^1$. This surface was constructed as $ \mathbb{Z}_{2}^4$-cover of $ \mathbb{P}^1\times \mathbb{P}^1$ (see \cite{MR4008073}). It is worth pointing out that C. Gleissner, R. Pignatelli and C. Rito constructed a family of surfaces with $ K_X^2 = 24, p_g\left( X\right) =3, q\left( X\right) =1 $ and $ d = 24 $ (\cite{2018arXiv180711854G}). Their example has a very similar construction as $ \mathbb{Z}_{2}^3$-cover of $ \mathbb{P}^1 \times C $ branched on ``fibers'' of the obvious trivial fibrations.    
    
    Throughout this note all surfaces are projective algebraic over the complex numbers. The linear equivalence of divisors is denoted by $ \equiv $. A character $ \chi $ of the group $ G$ is a homomorphism from $ G$ to $ \mathbb{C}^{*} $, the multiplicative group of the non-zero complex numbers. The rest of the notation is standard in algebraic geometry.
    	
    \section{$ \mathbb{Z}_{2}^3$-coverings}\label{z_2^3-coverings}
    The construction of abelian covers was studied by R. Pardini in \cite{MR1103912}. For details about the building data of abelian covers and their notations, we refer the reader to Section 1 and Section 2 of R. Pardini's work (\cite{MR1103912}). For the sake of completeness, we recall some facts on $ \mathbb{Z}_{2}^3 $-covers, in a form which is convenient for our later constructions.
 
   We denote by  $ \chi_{j_1j_2j_3} $ the character of $ \mathbb{Z}_{2}^3 $ defined by
   \begin{align*}
   	\chi_{j_1j_2j_3}\left( a_1,a_2,a_3\right): =  e^{\left( \pi a_1j_1\right) \sqrt{-1}}e^{\left( \pi a_2j_2\right) \sqrt{-1}}e^{\left( \pi a_3j_3\right) \sqrt{-1}}
   \end{align*}
   for all $ j_1,j_2,j_3,a_1,a_2,a_3\in \mathbb{Z}_2 $. A $ \mathbb{Z}_{2}^3 $-cover $ \xymatrix{X \ar[r] & Y} $ can be determined by a collection of non-trivial divisors $ L_{\chi} $ labelled by characters of $ \mathbb{Z}_{2}^3 $ and effective divisors $ D_{\sigma} $ labelled by non-trivial elements of $ \mathbb{Z}_{2}^3 $ of the surface $ Y $. More precisely, from \cite[\rm Theorem 2.1]{MR1103912} we can define $ \mathbb{Z}_{2}^3 $-covers as follows:
   \begin{Proposition} \label{Construction of cover of degree 8}
   	Given $ Y $ a smooth projective surface, let $ L_{\chi} $ be divisors of $ Y $ such that $ L_{\chi} \not\equiv \mathcal{O}_Y $ for all non-trivial characters $ \chi $ of $ \mathbb{Z}_{2}^3  $ and let $ D_{\sigma} $ be effective divisors of  $ Y $ for all $ \sigma \in \mathbb{Z}_{2}^3 \setminus \left\lbrace \left(0,0,0 \right)  \right\rbrace  $ such that the total branch divisor $ B:=\sum\limits_{\sigma \ne 0}{D_{\sigma}} $ is reduced. Then $ \left\lbrace L_{\chi}, D_{\sigma} \right\rbrace_{\chi,\sigma}$ is the building data of a $ \mathbb{Z}_{2}^3$-cover $ \xymatrix{f:X \ar[r]& Y} $ if and only if
   	\begin{align}\label{The condition of Z_2^3 covers} 
   		L_{\chi}+L_{\chi'} \equiv L_{\chi\chi'}+ \sum\limits_{\chi\left( \sigma\right)=\chi'\left( \sigma\right) = -1 }{D_{\sigma}}	
   	\end{align}
   	for all non-trivial characters $ \chi, \chi' $ of $ \mathbb{Z}_{2}^3  $. 
   	
   \end{Proposition}
   \noindent
   For the reader's convenience, we leave here the relations (\ref{The condition of Z_2^3 covers}) of the reduced building data of $ \mathbb{Z}_{2}^3 $-covers:
   $$
   \begin{adjustbox}{max width=\textwidth}
   	\begin{tabular}{l l l r r r r r r}
   		$ L_{100}+L_{100} $&$ \equiv $&$  $&$ $&$ D_{100} $&$ +D_{101 } $&$ +D_{110} $&$ +D_{111} $& \\
   		$ L_{100}+L_{010} $&$ \equiv $&$  $&$  $&$  $&$ $&$ D_{110} $&$ +D_{111} $& $ +L_{110} $\\
   		$ L_{100}+L_{001} $&$ \equiv $&$ $&$  $&$$&$ D_{101} $&$  $&$ +D_{111 } $& $ +L_{101} $\\
   		$ L_{010}+L_{010} $&$ \equiv $&$ D_{010 } $&$ +D_{011} $&$  $&$  $&$ D_{110} $&$ +D_{111 } $&\\
   		$ L_{010}+L_{001} $&$ \equiv $&$ $&$ D_{011} $&$  $&$  $&$  $&$ +D_{111 } $& $ +L_{011} $\\
   		$ L_{001}+L_{001} $&$ \equiv D_{001} $&$  $&$ +D_{011} $&$  $&$ +D_{101} $&$  $&$ +D_{111} $.&
   	\end{tabular}
   \end{adjustbox}
   $$	
   
   \noindent
   By \cite[\rm Theorem 3.1]{MR1103912} if each branch component $D_\sigma$ is smooth and the total branch locus $B $ is a simple normal crossings divisor, the surface $X$ is smooth. \\
   
   \noindent
   Also from \cite[\rm Lemma 4.2, Proposition 4.2]{MR1103912} we have:
   \begin{Proposition}\label{invariants of Z_2^3 cover}
   	If $ Y $ is a smooth surface and $ \xymatrix{f: X \ar[r]& Y} $ is a smooth $  \mathbb{Z}_{2}^3$-cover with the building data $ \left\lbrace L_{\chi}, D_{\sigma} \right\rbrace_{\chi,\sigma}$, the surface $ X $ satisfies the following:
   	\begin{align*}
   		2K_X & \equiv f^*\left( 2K_Y + \sum\limits_{\sigma \ne 0} {D_{\sigma} } \right); \\
   		f_{*}\mathcal{O}_X &= \mathcal{O}_Y \oplus \bigoplus\limits_{\chi \ne \chi_{000}  }L_{\chi}^{-1}.
   	\end{align*}	
   	\noindent
   	This implies that	
   	\begin{align*}
   		H^{0}\left( X, K_X\right) &= H^{0}\left( Y, K_{Y}\right) \oplus \bigoplus_{\chi \ne  \chi_{000}}{H^{0}\left( Y, K_{Y} +L_{\chi}\right)};\\
   		K^2_X &= 2\left( 2K_Y + \sum\limits_{\sigma \ne 0} {D_{\sigma} } \right)^2; \\
   		p_g\left( X \right) &=p_g\left( Y \right) +\sum\limits_{\chi \ne  \chi_{000}  }{h^0\left( L_{\chi} + K_Y \right)}; \\
   		\chi\left( \mathcal{O}_X \right) &= 8\chi\left( \mathcal{O}_Y \right)  +\sum\limits_{\chi \ne \chi_{000}  }{\frac{1}{2}L_{\chi}\left( L_{\chi}+K_Y\right)}. 
   	\end{align*}
   \noindent
   Moreover, the canonical linear system $ \left|  K_X \right|  $ is generated by
   \begin{align*}
   	f^{*}\left(  \left| K_Y + L_{\chi}\right| \right) +\sum\limits_{\chi\left( \sigma\right)=1 }{R_{\sigma}}    , \hskip 5pt \forall \chi \in J 
   \end{align*}
   \noindent
   where $ J:= \left\lbrace  \chi' : \left| K_Y + L_{\chi'}\right| \ne \emptyset \right\rbrace  $ and $ R_{\sigma} $ is the reduced divisor supported on $ f^{*}\left( D_{\sigma} \right)  $.
   \end{Proposition}
   \noindent
   For the proof of the last statement of Proposition \ref{invariants of Z_2^3 cover}, we refer the reader to \cite[\rm Page 3]{2018arXiv180711854G}.	
   
    \section{Construction}
    Throughout this section, we denote by $ Y:= \mathbb{P}^1 \times C $ the product surface of the rational curve $ \mathbb{P}^1 $ and a smooth elliptic curve $ C $. Let $ \xymatrix{p_1: \mathbb{P}^1 \times C \ar[r] &\mathbb{P}^1} $ be the projection of the product surface $ \mathbb{P}^1 \times C $ on $ \mathbb{P}^1 $. We denote by $ E $ a general fiber of $ p_1 $. The canonical class of $ Y $ is $ K_Y \equiv -2E $.\\
    
    Let $ E_1, E_2, \ldots, E_6 $ be distinct elliptic fibres and let $ F_{1},F_{2}, \ldots, F_{n}$, $ F_{1}^{'}$, $ F_{2}^{'}, \ldots, F_{n}^{'} $, $ F_{1}^{''}, F_{2}^{''}$, $ F_{3}^{''}$ be distinct rational fibres (with $ n \ge 3 $) such that $ 2F_{1}^{''} \equiv 2F_{2}^{''} \equiv 2F_{3}^{''} $. Because the sum of two points in an elliptic curve is divisible by $ 2 $ in the Picard group, there are fibres $ F_{ii} $ such that $ 2F_{ii} \equiv F_{i} + F_{i}^{'} $, for all $ i\in \left\lbrace 1,2,\ldots,n \right\rbrace  $. We consider the following divisors	
    $$
    \begin{adjustbox}{max width=\textwidth}
    	\begin{tabular}{l l l l}
    		$ D_{100}: = E_1 + E_2$, & $ D_{101}: = E_3 + E_4$, & $ D_{110}: = E_5 + E_6$, & $ D_{111}: = \sum\limits_{i=1}^{n}{\left( F_{i}+F_{i}^{'}\right) } $ \\
    		$ L_{100}:=3E + \sum\limits_{i=1}^{n}{F_{ii}} $&$ L_{010}:=E + \sum\limits_{i=1}^{n}{F_{ii}} +\eta_1 $&$ L_{001}:=E + \sum\limits_{i=1}^{n}{F_{ii}} +\eta_2$& $ L_{110}:=2E +\eta_1 $ \\
    		$ L_{101}:=2E +\eta_2 $ &$ L_{011}:=2E +\eta_3 $ &$ L_{111}:=E + \sum\limits_{i=1}^{n}{F_{ii}} +\eta_3$,&
    	\end{tabular}
    \end{adjustbox}
    $$
    \noindent
    where $ \eta_1:= F_{1}^{''} - F_{2}^{''}$, $ \eta_2:= F_{2}^{''} - F_{3}^{''}$,  $ \eta_3:= F_{1}^{''} - F_{3}^{''} = \eta_1 + \eta_2$ are non-trivial $ 2 $-torsions. These divisors $ D_{\sigma}, L_{\chi} $ satisfy the following relations
    $$
    \begin{adjustbox}{max width=\textwidth}
    	\begin{tabular}{l l r r r r r r l}
    		$ L_{100}+L_{100} $&$ \equiv $&$ D_{100} $&$ +D_{101 } $&$ +D_{110} $&$ +D_{111} $& &$ \equiv $&$ 6E + \sum\limits_{i=1}^{n}{2F_{ii}} $\\
    		$ L_{100}+L_{010} $&$ \equiv $&$  $&$ $&$ D_{110} $&$ +D_{111} $& $ +L_{110} $&$ \equiv $&$ 4E + \sum\limits_{i=1}^{n}{2F_{ii}} +\eta_1 $\\
    		$ L_{100}+L_{001} $&$ \equiv $&$$&$ D_{101} $&$  $&$ +D_{111 } $& $ +L_{101} $&$ \equiv $&$ 4E + \sum\limits_{i=1}^{n}{2F_{ii}} +\eta_2 $\\
    		$ L_{010}+L_{010} $&$ \equiv $&$  $&$  $&$ D_{110} $&$ +D_{111 } $&&$ \equiv $&$ 2E + \sum\limits_{i=1}^{n}{2F_{ii}} $\\
    		$ L_{010}+L_{001} $&$ \equiv $&$  $&$  $&$  $&$ +D_{111 } $& $ +L_{011} $&$ \equiv $&$ 2E + \sum\limits_{i=1}^{n}{2F_{ii}} +\eta_3 $\\
    		$ L_{001}+L_{001} $&$ \equiv $&$  $&$ D_{101} $&$  $&$ +D_{111} $&&$ \equiv $&$ 2E + \sum\limits_{i=1}^{n}{2F_{ii}} $.
    	\end{tabular}
    \end{adjustbox}
    $$	
    \noindent
    Thus by Proposition \ref{Construction of cover of degree 8}, the divisors $ D_{\sigma}, L_{\chi} $ define a $ \mathbb{Z}^3_2 $-cover $ \xymatrix{f: X \ar[r] & Y}  $. Because each branch component $ D_{\sigma} $ is smooth and the total branch locus $ B $ is a normal crossings divisor, the surface $ X $ is smooth. Moreover, by Proposition \ref{invariants of Z_2^3 cover}, the surface $ X $ satisfies the following
    \begin{align*}
    	2K_X & \equiv f^*\left( 2E + \sum\limits_{i=1}^{n}{\left( F_{i}+F_{i}^{'}\right) }\right).
    \end{align*}
    \noindent
    We notice that a surface is of general type and minimal if the canonical divisor is big and nef (see e.g. \cite[\rm Section 2]{MR2931875}). Since the divisor $ 2K_X $ is the pull-back of a nef and big divisor, the canonical divisor $ K_X $ is nef and big. Thus, the surface $ X $ is of general type and minimal. Moreover, by Proposition \ref{invariants of Z_2^3 cover} the invariants of $ X $ are as follows
    \begin{align*}
    	K_X^2= 16n, p_g\left( X\right)= 2n, \chi\left( \mathcal{O}_X\right) =2n, q\left( X\right)=1.
    \end{align*}
    
    We show that the canonical map is of degree $ 8 $. By Proposition \ref{invariants of Z_2^3 cover}, we have the following decomposition     
    \begin{align*}
    	H^{0}\left( X, K_X\right) = H^{0}\left( Y, K_{Y}\right) \oplus \bigoplus_{\chi \ne  \chi_{000}}{H^{0}\left( Y, K_{Y} +L_{\chi}\right)}. 
    \end{align*}
    \noindent
    Moreover, the choice of $ L_{\chi} $'s yields that
    \begin{align*}
    	h^0\left( L_{ \chi} + K_Y \right)  = 0
    \end{align*}
    for all $ \chi \ne \chi_{100}  $. By Proposition \ref{invariants of Z_2^3 cover}, the linear system $ \left| K_X \right| $ is generated by
    \begin{align*}
    	f^{*}\left( K_Y + L_{100}\right)  = f^{*}\left( E + \sum\limits_{i=1}^{n}{F_{ii}}\right).
    \end{align*}
    \noindent
    This implies that the canonical map of the surface $ X $ factors through the $ \mathbb{Z}^3_2 $-cover $ \xymatrix{f: X \ar[r] & Y}$. Thus the following diagram commutes
    $$
    \xymatrix{X \ar[0,2]^{\mathbb{Z}_2^3}_f \ar[2,2]_{\varphi_{\left| K_X \right| }}&& Y \ar[2,0]^{\varphi_{\left| E + \sum\limits_{i=1}^{n}{F_{ii}} \right| }}\\
    &&\\
    &&\mathbb{P}^{2n-1}.}
    $$
    \noindent
    Since the map $ \varphi_{\left| E + \sum\limits_{i=1}^{n}{F_{ii}} \right| } $ is of degree one for all $ n \ge 3 $, the canonical map of $ X $ is of degree $ 8 $.\\
    
    \begin{Remark}
    	Let $ Z:= X/\Gamma $ be the quotient surface of $ X $, where $ \Gamma:= \left\langle \left( 0,0,1\right), \left( 0,1,0\right) \right\rangle  $ is the subgroup of $ \mathbb{Z}_2^3 $. The surface $ Z  $ is a surface of general type whose only singularities are $ 12n $ nodes. Moreover, the canonical map $ \varphi_{\left| K_{Z} \right| } $ is a map of degree $ 2 $ (see \cite[\rm Theorem 5.1]{MR656051} and \cite[\rm Theorem 1.1]{MR3272911}). The canonical map $ \varphi_{\left| K_X \right| } $ is the composition of the degree $ 4 $ quotient map $ \xymatrix{X \ar[r]& Z:= X/\Gamma } $ with the canonical map $\varphi_{\left| K_{Z} \right| }$ of $ Z $ (see e.g. \cite[\rm Example 2.1]{MR1103913}).
    \end{Remark}
    
    \begin{Remark}\label{deformationvariation}
    	In the above construction, there are four different possible choices for each $ F_{ii} $. A different choice produces a different surface $ X $. 
    \end{Remark}
    
    \noindent
    The previous observation leads to the following interesting question:
    \begin{Question}
    	Are the surfaces in Remark \ref{deformationvariation} deformation equivalent?
    \end{Question}

\section*{Acknowledgments}
The author is deeply indebted to Margarida Mendes Lopes for all her help. Thanks are also due to Jungkai Alfred Chen for many interesting conversations. The author would like to express his gratitude to the anonymous referee for his/her thorough reading of the paper and valuable suggestion. This paper was finished during the author's postdoctoral fellowship at the National Center for Theoretical Sciences (NCTS), Taiwan, under the grant MOST 110-2123-M-002-005. The author would like to thank NCTS for the financial support and kind hospitality.

%\nocite{*}
%\bibliographystyle{acm}
%\bibliography{nguyenbin}

\Addresses

\end{document}